# ISOMETRIES GROUPS AND A MULTIRESOLUTION ANALYSIS ON SUB-RIEMANNIAN MANIFOLDS

Romina Cardo    &    Alvaro Corvalán

**ABSTRACT.** In this letter we exhibit the relation between the isometries of a Riemannian contraction of a sub-Riemannian manifold and those of the sub-Riemannian metric, for to use this relation with two goals: establishing a result about the existence of fixed points of isometries groups; and the other, defining a Multiresolution Analysis (MRA) on sub-Riemannian manifolds that it will permit to obtain Haar's bases on the manifolds before mentioned.

## 1. Introduction

The Sub-Riemannian geometry is an area has been studied with high intensity at the last two decades.In certain form, it is a generalization of the Riemannian geometry, that arose on having considered variational questions, where there were "prohibited" directions for which the length or the energy were infinite. Lately, the dynamical systems, the theory of control and the problems of pursuit and evasion; they motivate the work in sub-Riemannian geometry.In this letter, we study the existence of fixed points of certain groups of isometries for sub-Riemannian manifolds and we introduce a Multiresolution Analysis (MRA) on sub-Riemannian manifolds.We start by giving some definitions and preliminary results:

## 2. Sub-Riemannian Manifolds

**Definition 2.1.** A *Sub-Riemannian Manifold* $(M, g)$ consists of a differential manifold $M$ with $n = \dim(M)$ (we will suppose $n \geq 3$) and an application $g : TM^* \to TM$ smooth such that $g_p : M_p^* \to M_p$ is linear, symmetric and non-negative definite,that is:
   (♦) If $\theta \in M_p^* \Rightarrow \theta(g_p(\theta)) \geq 0$ with $g_p = g / M_p$
   (♦♦) If $\theta, \xi \in M_p^* \Rightarrow \theta(g_p(\xi)) = \xi(g_p(\theta))$
   provided with a vector subbundle $S$ of the tangent bundle $TM$, whose fibers have fixed positive codimension $n - m$, and such that $TM$ is the vector subbundle that it contains $S$ and all the tangent vectors obtained by the restrictions in each fiber of the Lie brackets of the sections of $S$; with a differential application $Q$ from $S \times S$ on the real nonnegative numbers, that it is bilinear and definite positive, named *Sub-Riemannian metric*.
   We suppose that $S$ verifies the *Strong bracket generating hypothesis;* we mean by this, that for all $p \in M$ it fulfills the *2- Hörmander´s Condition* for each $p$, that is: for any $v \in S_p$ (the fiber over $p$) non-null, such that: $S_p + [v, S_p] = M_p$, it is to say, if $S = \mathrm{Im}(g) \subseteq TM$ with $r = \dim(S)$, $S$ is a subbundle of constant codimension $n - r$ and the following flag of subspaces of $M_p$ finishes at $M_p$:

$S_p \subseteq S_p^2 \subseteq S_p^3 \subseteq \ldots \subseteq M_p$ where
$S_p^2 := S_p + [S_p, S_p]$
$S_p^3 := S_p^2 + [S_p, S_p^2]$
$S_p^4 := S_p^3 + [S_p, S_p^3]$, etc.
are the subspaces spanned from $S_p^l$ and the Lie brackets of the sections of $S_p$ and $S_p^l$.

Most questions, and this one too, are easier if $S_p^2 = M_p$ (Strong Bracket Generating Hypothesis-the2 Hörmander Condition).

This induces a smooth varying, positive definite, bilinear form on $S_p$:
$\langle g_p(\theta), g_p(\xi) \rangle = \theta(g_p(\xi))$ (if $n - r = 0$, we have the Riemannian case) and taking basis of $M_p^*$ and $M_p$ then $\left(g_p^{ij}\right)_{1 \leq i,j \leq n}$ is the symmetric matrix (ie, $g_p^{ij} = g_p^{ji}$) of the raised indices like the Riemannian case.

If we take a complementary subbundle $S^\perp$ and make it orthogonal to $S$, we have got a Riemannian metric, and we talk of a "Riemannian Contraction of the Sub-Riemannian Structure", immediately we will see the reason of the word "contraction".

If we have the Strong Bracket Generating Hypothesis and we take $\{X_i\}$ an orthonormal moving "frame" that spans $S$ and $S^\perp = \mathrm{span}_{i_h, i_j \in I} \{[X_{i_h}, X_{i_j}]\}$; and we obtain the Riemannian contraction making $\{[X_{i_h}, X_{i_j}]\}$ orthonormal, then we can show that the measure given from the Riemannian contraction to the Sub-Riemannian manifold doesn't depends on $\{X_i\}$.

If $S$ is spanned by brackets, "a Sub-Riemannian metric $S$" it is a bilinear symmetric form $Q_p : S_p \times S_p \to \mathbf{R}$ in both coordinates differentiable at $p$, and positive definite over $S_p$. The sub-Riemannian metric $Q_p$ can be given too by means of the operator $g_p$ that we will define soon:

Using the Riesz´s Representation Theorem for Bundles, we can define the linear functional $g_p : M_p^* \to M_p$ (with image $S_p$), and its associated bilinear form $\phi_p$, that results symmetric and nonnegative definite, such that $Q_p(v, g_p(\varphi)) = \varphi(v)$; $\forall v \in S_p$.

Now, we will assume that the sub-Riemannian metric is given by means of $g_p$.

**Definition 2.2.** We named $g$ to the *sub-Riemannian metric of $M$*, with $g_p := g \mid M_p$; then $g$ provides the manifold with a structure of metric space (which distance we will write $d$) compatible with its own topology by the CHOW´s Theorem [5].

**Definition 2.3.** We name *Riemannian Contraction $f$* of the sub-Riemannian metric $g$, to a Riemannian metric for $M$ which restriction to $S$ realize the sub-Riemannian metric $g$ of $M$.

### 3. Preliminary Foundations

**Definition 3.1.** We say that $\Psi : M \to M$ is an *isometry* if $\Psi$ preserves distances, that is, $d(\Psi(x); \Psi(y)) = d(x, y)$; $\forall x, y \in M$; and we say that $\Psi$ is $C^2$-*infinitesimal* if $\Psi$ is $C^2$, and, $g(\Psi(x)) = d\Psi(x) \cdot g(x) \cdot d\Psi^*(x)$. It can be proven ( Strichartz [12] ) that in these conditions, $\Psi$ commutes with the exponential mapping, that is to say :
$\Psi\left(\exp_p(u)\right) = \exp_{\Psi(p)}(d\Psi(p)^* u)$ ; $\forall p \in M$, and, $\forall u \in M_p^*$.

**3.2.** There must happen that: $[d\Psi] \cdot [g] \cdot [d\Psi^*] = [g]$, with $[d\Psi^*] = [d\Psi]^t$, where $[\ ]$ denotes the matrices of the respective linear transformations.

In the same way, being $\Psi$ an infinitesimal isometry for the contraction $f$, we have :
$[d\Psi] \cdot [f] \cdot [d\Psi^*] = [f]$.

In the Riemannian case, we can choose a basis $B$ so that $[f]_B = Id = \begin{bmatrix} I & 0 \\ 0 & I \end{bmatrix}$.

For the sub-Riemannian case we will have that $[g]_B = \begin{bmatrix} I & 0 \\ 0 & 0 \end{bmatrix}$, where $range([g]_B) = \dim(S)$.

So we will obtain that $[d\Psi] \cdot \begin{bmatrix} I & 0 \\ 0 & I \end{bmatrix} \cdot [d\Psi]^t = \begin{bmatrix} I & 0 \\ 0 & I \end{bmatrix}$ then $[d\Psi] \cdot \begin{bmatrix} I & 0 \\ 0 & 0 \end{bmatrix} \cdot [d\Psi]^t = \begin{bmatrix} I & 0 \\ 0 & 0 \end{bmatrix}$ it is to say, $[d\Psi] \cdot [g] \cdot [d\Psi]^t = [g] \quad \Rightarrow \quad d\Psi \circ g \circ d\Psi^* = g$ .

Hence, if $\Psi$ is an infinitesimal isometry for $f$, then it is an isometry for $g$ too.

### 3.3. THEOREM of ESTIMATION of METRICS ( Strichartz [12] )

If $S$ verifies the 2-Condition of Hörmander and $d_R$ is the associated metric to the Riemannian contraction $f$, then it exists a positive constant $c$ such that: $d(x,y) \leq c \cdot (d_R(x,y))^{\frac{1}{2}}$ ; $\forall x, y \in K$ where $K$ is any compact of $M$.

### 3.4. COROLLARY of THE PREVIOUS THEOREM:

Let $\Psi : M \to M$ an infinitesimal isometry for the Riemannian contraction $f$ ( $\therefore$ it is an infinitesimal isometry for $g$ from 3.2.) and if $p$ is a fixed point for $f$, then it is an infinitesimal isometry for $g$ too.

**Definition 3.5.** Let $p \in M$, we named *radius of Strongly Convergence for $p$*, that we will write $\xi_p$, so:
$\xi_p := \sup\{r \in \mathbf{R}_{\geq 0} \cup \{+\infty\}$ such that $\forall\ q, s \in B_r(p)$ ,$\exists!$ minimizing geodesic that joins them contained in $B_r(p)\}$

**3.6.** We conclude with the principal result:

Let $p \in M$, for a known consequence of *Riemann's Principle [8]*, we can affirm the existence of a strongly convex neighbourhood for $f$, then $\xi_p > 0$ (for each $p$ ).

Let $H$ a compact group of isometries for $f$ ( $\therefore$ it is for $g$ for the above mentioned fact previously in 3.2. ) such that $diam(H.p) \leq \frac{1}{2}.\xi_p$ ( and $\therefore$  $H.p \subseteq B_{2.(\frac{1}{2}.\xi_p)}(p) = B_{\xi_p}(p)$ ).

Then we obtain that the orbit $H.p \subset B_{\frac{1}{2}.\xi_p}(p) \subset B_{\xi_p}(p)$ is strongly convex. Using the *KLINGENBERG-EBERLEIN's THEOREM [8]*: $H$ have a fixed point in $B_{\frac{1}{2}.\xi_p}(p)$ for $f$, that for the previous Corollary 3.4, it is fixed point for $g$ too , then $H$ have a fixed point for $g$ .

**NOTE.** In particular: If the manifold is compact, it is enough to ask that $H$ is closed.

**3.7.** In the previous conditions, as the isometries of $H$ are $C^2$ then they are regular, so now our result allows us to affirm that $H$ results isomorphic to a compact subgroup of group $O(m)$ of orthogonal matrices of dimension $m = \dim(S)$ .(this is a consequence of Strichartz [12] 8.5 pág. 248 ).

    **3.8.** The key of the previous result is to pass to the Riemannian contraction, because for sub-Riemannian metrics never result the *Riemann's Principle* , it is to say, strongly convex neighbourhoods of a point do not exist, because, all point $q$ of a sub-Riemannian manifold is an accumulation point of its "CUT LOCUS" , being this one the set of all such points that there exists more than one minimizing geodesic that joins them with $q$.(Strichartz [12] page 260, theorem 11.3 and its previous paragraph).

    **3.9.** The result of this letter arose of trying to prove the analogous one to the *KLINGENBERG-EBERLEIN'S THEOREM [8]* for sub-Riemannian manifolds.

## 4. Lattices on Sub-Riemannian Manifolds

    **Definition 4.1.** We say that $K$ is a **Lattice for M** if $K$ is a countable subgroup of isometries of $M$.
    **Definition 4.2.** We say that $K$ is an **Uniform Lattice for M** if $K$ is a lattice for $M$ and $M/\sim$ is compact where $\sim$ is the equivalence relation: $x \sim y$ if and only if $\exists J \in K : y = J(x)$.
    **Definition 4.3.** Let **A** a diffeomorphism from **M** to **M**, we say that $K$ is an **Invariant Lattice with respect to A** if $K$ is a lattice for $M$ and $\forall J \in K \ \exists \tilde{J} \in K : A \circ J \circ A^{-1} = \tilde{J}$ (that is $A \circ J = \tilde{J} \circ A$ we denote $a(J) = A \circ J \circ A^{-1}$ that results an endomorphism of $K$ ).
    **Definition 4.4.** Let $M$ be an oriented Sub-Riemannian or Riemannian Manifold with $\dim(M) \geq 3$, $K$ a group of isometries of $M$ , this induces an action: $K \times M \to M$ defined by: $(\tau,p) \mapsto \tau(p)$, $A : M \to M$ a diffeomorphism such that:
    $(a) \ \forall J \in K, \exists \tilde{J} \in K / A \circ J = \tilde{J} \circ A$; so we have a group homomorphism
        $a : K \to K$ with $a(J) = A \circ J \circ A^{-1}$ .
    $(b)$ with the norm induced in the fibers of $TM$, the tangent bundle, from the Riemannian structure (or a Riemannian contraction in the sub-Riemannian case):
$$\sup_{p \in M} \|dA^{-1}\|_{M_p} \leq c < 1 ,$$
For then we obtain:
$$d(A^{-1}(p), A^{-1}(q)) < c \cdot d(p,q)$$
We will say that $A$ is a **dilation** and from $A$ we obtain the scaling process.

## 5. A Multiresolution Analysis (MRA) on Sub-Riemannian Manifolds

    Let $M$ be an oriented manifold endowed with a Sub-Riemannian metric and a measure $dm$ :
    There are different ways to choose $dm$, actually, as a sub-Riemannian manifold has a metric space structure, so there are Hausdorff measures.
    Notwithstanding, we will take $dm$ the Riemannian Volume (or the Sub-Riemannian one,

obtained from a certain Riemannian contraction of the Sub-Riemannian structure in a way we will see later).

**Definition 5.1.** In the context of being $K$ an uniform lattice for $M$ and $A$ a dilation, a **MULTIRESOLUTION ANALYSIS** (*MRA*) is:

A sequence $\{V_i\}_{i \in \mathbf{Z}}$ of closed subspaces of $L^2(M, dm)$ such that:

(✠) $V_i \subset V_{i+1}$ , $\forall i \in \mathbf{Z}$ .

(✠) $f \in V_i \Leftrightarrow f \circ A^{-i} \in V_0$ ; $f \in V_0 \Leftrightarrow f \circ J \in V_0$ , $\forall J \in K$ .

(✠) $\overline{\underset{i \in \mathbf{Z}}{\cup} V_i} = L^2(M, dm)$ .

(✠) $\underset{i \in \mathbf{Z}}{\cap} V_i = \{0\}$ .

(✠) $\exists \varphi \in V_0$ (the generator of $\{V_i\}_{i \in \mathbf{Z}}$) such that:

(♠) $V_0 = \overline{span \{\varphi \circ J , J \in K\}}$

(♠♠) $\exists\; \alpha_1, \alpha_2 > 0$ such that if $\{\lambda_j\}_{j \in K} \in l_2(K)$ then

$$\alpha_1 . \|\lambda\|^2_{l_2(K)} \leq \|\sum_{j \in K} \lambda_j . (\varphi \circ J)\|^2_{l_2(K)} \leq \alpha_2 . \|\lambda\|^2_{l_2(K)}.$$

**5.2.** Isometries and dilations we need a metric space structure.

We defined a Sub-Riemannian Metric, the construction is as follows:

We say that $\gamma : [a,b] \to M$ is a **Lengthy** curve if $\gamma$ is continous, piecewise $C^1$ and if $\exists \xi : [a,b] \to TM^*$ such that $\xi(t) \in M^*_{\gamma(t)}$ and $g_{\gamma(t)}(\xi(t)) = \dot{\gamma}(t)$ that we say a **cotangent lift**.

We observe that, then the direction of $\gamma$ belongs to the **admissible** bundle $S$ at every point. Then:

the length of $\gamma$ is $l(\gamma) = \int_a^b \left( \xi(t)\left(\dot{\gamma}(t)\right) \right)^{\frac{1}{2}} dt = \int_a^b (\xi(t)(g_{\gamma(t)}(\xi(t))))^{\frac{1}{2}} dt$ .

(▲) From the Theorem of CHOW [5], for $p$ and $q$ any points of $M$, there is a lengthy curve joining $p$ and $q$ .

The key idea comes from the fact that $e^{-tY}.e^{-tX}.e^{tY}.e^{tX} = e^{t^2[X,Y]+o(t^3)}$ for $X$ and $Y$ sections of $TM$, and then, if $\text{Im} X \subseteq S$, $\text{Im} Y \subseteq S$, and $S_p \subseteq S_p^2 \subseteq S_p^3 \subseteq \ldots \subseteq M_p$, we can travel in the lost directions turning and going up.

(▲) We define the distance between $p$ and $q$ points of $M$ by:

$d(p,q) := \inf \{l(\gamma) = $ length of $\gamma$, with $\gamma$ any lengthy curves joining $p$ and $q \} < +\infty$

(▲) Analogous of many theorems of Riemannian Manifolds [2] can be obtained as the Gauss Lemma, Hopf-Rinow theorem, solutions of the Hamilton-Jacobi equations are locally minimizing curves named regular goedesics.

(▲) The reason of the name: **Riemannian Contraction** is clear now, because there are more lengthy curves between $p$ and $q$, because if $\gamma$ is a lengthy Sub-Riemannian curve, $l(\gamma) = l_R(\gamma)$ : the Riemannian length so $d_R(p,q) \leq d(p,q)$ .

(▲) An **isometry** $\Psi : M \to M$ is a preserving distance application, ie, $d(\Psi(p), \Psi(q)) = d(p,q)$

(▲) Robert Strichartz proved in 1985 that if $\Psi \in C^2$ , $\Psi$ is an isometry and if $g_{\Psi(p)} = d\Psi(p).g_p.d\Psi^*(p)$ then $\Psi$ commute with the exponential application.

(▲) If $\Psi : M \to M$ is an isometry for a Riemannian contraction as the one mentioned before (from a moving frame of $S$), we found that it's an isometry for the Sub-Riemannian structure too (and it is $C^1$ by the Myers-Steenrod theorem) .

## 6. GOAL

(∗) If $(M,g)$ is a Sub-Riemannian Manifold, we study conditions that allow us to give a **Fundamental and Self-similar Set** with positive measure respect of a group of isometries and a dilation that it is a conjugation against the group of isometries. With this, it is possible to define an **MULTIRESOLUTION ANALYSIS** (*MRA*) for $L^2(M)$

(∗) The idea is to take advantage partly of the construction done by Stephan Dahlke [6] for Riemannians Manifolds.

(∗) Let $K$ a group of infinitesimal isometries of $(M,g)$, that for a result of Robert Strichartz [12], they are also isometries in the sense of preserving the induced distance.

(∗) In particular, if $(M,f)$ is a Riemannian contraction of $(M,g)$, the isometries of $(M,f)$ are isometries of $(M,g)$; so that we can use the groups of isometries of the Riemannian contractions.

(∗) An important part of the construction that follow can be done in metric complete spaces and locally compact in general. There are two advantages on working with manifolds:

(1) It is possible to use the Theorem of Change of Variables [14] to prove that the characteristic function of a self-similar set generates a **MULTIRESOLUTION ANALYSIS** (*MRA*) for $L^2(M)$.

(2) There are more comfortable ways of finding isometries and of working with them.

## 7. CONDITIONS

Let $K$ a countable group of isometries of $(M,g)$ a complete sub-Riemannian manifold, such that:

(*i*) $\forall p \in M \; \exists U_p$ neighbourhood of $p$ such that $J(U_p) \cap U_p = \emptyset$ for every except finite $J \in K$.

(*ii*) If $p$ and $q$ do not belong to the same orbit ($J(p) \neq q, \forall J \in K$) $\exists U_p, U_q$ neighbourhoods of $p$ and $q$, respectively, such that $J(U_p) \cap U_q = \emptyset, \forall J \in K$.

(*iii*) The set of points of $M$ that there are fixed for any $J \in K$ ($J \neq Id$) is never dense in $M$.

Additionally we will require:

(♣) $M/\sim$ is compact where $\sim$ is the equivalence relation which classes are the orbits of the elements provided by $K$.

The following construction is essentially a translation from the one given by Stephan Dahlke [6] for the Riemannian case.

The idea, here, is to see that sometimes we can obtain results for the Sub-Riemannian case from the Riemannian case.

**7.1. Important Observation.** Stephan Dahlke [6], in his work, mentioned to "the fixed points of $K$", that usually means "fixed for all elements of $K$". But if it exists $p_0$ / $J(p_0) = p_0, \forall J \in K$, then we will have a contradiction with (*i*) because $J(U_{p_0}) \cap U_{p_0} \supseteq \{p_0\} \neq \emptyset, \forall J \in K$. So we will use (*iii*) as we enunciated it.

## 8. FUNDAMENTAL SETS

**Definition 8.1.** A set $F \subseteq M$ is named **FUNDAMENTAL for** $K$ if verify:

(*a*) $F = \overline{F^\circ}$ (is the closure of its interior)

(*b*) $\bigcup_{J \in \mathbf{K}} J(F) = M$

(*c*) $J(F^\circ) \bigcap F^\circ = \emptyset$ , $\forall J \neq Id$

**Definition 8.2.** A set $F \subseteq M$ is called **SELF-SIMILAR** if there exists a finite set $J_i \in K$ such that $A(Q) = \bigcup_{i \in I} J_i(Q)$, where $A$ is the inverse of $A^{-1}$, with $A^{-1}$ a contraction, that is to say $d(A^{-1}(p), A^{-1}(q)) < c.d(p,q)$ $\forall p, q \in M$.

Stephan Dahlke [6] enunciates but does not demonstrate (though it is used) the existence of fundamental sets.

We will show how to construct them:

**Lemma 8.3.** We can choose first a $p_0 \in M$ that it is NOT fixed for any $J \in K$ with $J \neq Id$.

**Proof.** As $M$ is complete, $K$ is countable and the set remains fixed for each $J \neq Id$ is never dense, result that $M \neq \bigcup_{J \neq Id} \{\text{fixed points of } J\}$ because for the Baire's Theorem, a complete metric space is not union of countable never dense sets, then $\exists\, p_0 \in M$ that it is not fixed for any $J \neq Id$.

**Lemma 8.4.** Now, for $p_0$ and $\forall \epsilon > 0$ exists only finite images of $p_0$ for elements $J \in K$ with $d(p_0, J(p_0)) < \epsilon$.

**Proof.** We suppose that for some $\epsilon > 0$ it exists infinite (different) $J(p_0)$ in $B(p_0, \epsilon) \subseteq \overline{B(p_0, \epsilon)}$.

As we are on a manifold, and the differential manifolds, for the Whitney's Theorem [8], they can be imbedded in a $\mathbb{R}^n$ of dimension greater enough if the manifold is not compact, we obtain that $\overline{B(p_0, \epsilon)}$ is a closed bounded set in some euclidean space, and then, it is compact. If the manifold were compact, then $\overline{B(p_0, \epsilon)}$ is a closed set in a compact, and hence it is compact too.

So we would have a sequence $(J_n(p_0))_{n \in \mathbb{N}} \subseteq \overline{B(p_0, \epsilon)}$ (compact set), ( eventually taking a subsequence) that we can suppose that it has a limit $p \in \overline{B(p_0, \epsilon)}$.

But now, as $J_n(p_0) \xrightarrow[n \to +\infty]{} p$ results that in any neighbourhood $U_p$ of $p$ there are infinite different images of $p_0$.

Now if $p$ is not in the same orbit of $p_0$, we obtain a contradiction because we would have $U_{p_0}$ and $U_p$ such that $J(U_{p_0}) \cap U_p = \emptyset$, $\forall J \in K$, but $J_n(p_0) \in J_n(U_{p_0})$ and $J_n(p_0) \in U_p$.

And if $p_0$ and $p$ are in the same orbit $\exists J_0 \in K$ such that $J_0(p_0) = p$ but also we arrive to a contradiction because if we take a neighbourhood $U_p$ of $p$ such that $J(U_p) \cap U_p = \emptyset$ for almost all $J$, on the other hand we have that $(J_n \circ J_0^{-1})(U_p) \supseteq J_n(J_0^{-1}\{p_0\}) = J_n(\{p_0\} \subseteq U_p$, and hence, $(J_n \circ J_0^{-1})(U_p) \bigcap U_p \neq \emptyset$, for infinite $J = J_n \circ J_0^{-1}$. Then it can not have an $\epsilon > 0$ such that there are infinite images of $p_0$ with $d(p_0, J(p_0)) < \epsilon$.

## 9. Construction: POLYGONS of DIRICHLET

**Definition 9.1.** For the $p_0$ mencioned, and for each $J \in K$, $J \neq Id$, we define $H_q = \{p \in M \text{ such that } d(q,p_0)) < d(q,J(p_0))\}$ the open semispace respect of $J$, and the **Polygon of Dirichlet**:

$$D_{p_0} = \bigcap_{J \neq Id} H_J$$

and let $F = \overline{D_{p_0}}$. Then, $F$ is a fundamental set for $K$.

Coming back to the isometries and dilations, if we have a diffeomorphism $A : M \to M$ with $dA \subseteq S$, and, $\|dA^{-1}(p)\| \leq c < 1$ with the norm given from the Riemannian contraction of the Sub-Riemannian metric, we can use: the fact that both topologies coincides (this follows from any proof of Chow's theorem [5] ) and Stephan Dahlke's construction for the Riemannian manifolds [6]:

**9.2.** An algorithm for the construction of $Q$ a **fundamental and self-similar region**:

instead of the long procedure that proposes Dahlke [6] in his Lemma 3.1, it is sufficient to observe that $A^{-1} \circ J_i$ are contractions $\forall i \in I$ and by the **Iterated Function System (IFS) theorem** see [7]: Let $\{F_1, \ldots, F_m\}$ be a family of contractions on a complete metric space $X$ (that is an IFS) with $m \geq 2$, then there exists a unique, non-empty compact set $E \subset X$ that satisfies $E = \bigcup_{i=1}^{m} F_i(E)$, and then $E$ is the self-similar region that we need.

In an analogous way for the construction of a Multiresolution Analysis (MRA) we can follow literally the theorem 2.1 and the Lemma 3.2 of Dahlke [6]:

Supossing that exists a finite set $\{J_i\} \subset K$ such that $A(Q) = \bigcup_{i \in I} J_i(Q)$ where $A$ is the inverse of $A^{-1}$, and, $\|dA(p)\| \geq \alpha > 1$ and additionally requiring that $K/a(K)$ is finite with $a(J) = A \circ J \circ A^{-1}$.

Now, if there is a fundamental region $Q$ with measure $\mu(Q) > 0$, then taking $\varphi = \chi_Q$: the characteristic function of $Q$ - i.e. , the **Haar function**-, as the generator of $\{V_i\}_{i \in \mathbf{K}}$ and taking $V_0 = \overline{\text{span} \{\varphi \circ J, J \in K\}}$ then $\{V_i\}_{i \in \mathbf{K}}$ results a Multiresolution Analysis (MRA) for the Riemannian contraction and then, we get it for the Sub-Riemannian metric too.

Some tipical examples are Heisenberg type groups : $(\vec{x}, \vec{y})$ where $\vec{x}, \vec{y} \in R^n$ with the operation $(x,y) \circ (x',y') = (x + x', y + y' + L(x,x'))$ where $L$ is a skew-symmetric bilinear function.

These groups have got a natural Sub-Riemannian metric and homogeneous dilations $\delta_t(x,y) = (t.x, t^2.y)$, taking $A^{-1} = \delta_t$ with $0 < t < 1$. It's easy to verify that $A^{-1}$ satisfies the former hypothesis , and that the integer translations are a group of isometries both for the Riemannian contraction and the Sub-Riemannian metric, so we have a Multiresolution Analysis (MRA) for Sub-Riemannian Manifolds.

**Remark.** There is a related work of Robert Strichartz using a different approach on nilpotent Lie groups.